\documentclass{svproc}

\usepackage{url}

\usepackage{graphicx}              
\usepackage[caption=false]{subfig} 

\begin{document}
\mainmatter              
\title{A Two-Step Method Coupling Eddy Currents and Magneto-Statics}
\author{Martina Busetto
\and Christoph Winkelmann
}
%
%
\tocauthor{Martina Busetto and Christoph Winkelmann}
\institute{ABB Corporate Research, Segelhofstrasse 1K, 5405 Baden-D\"attwil, Switzerland
\email{martina.busetto@ch.abb.com, christoph.winkelmann@ch.abb.com}}

\maketitle              

\begin{abstract}
We present the mathematical theory and its numerical validation of a method tailored to include eddy-current effects only in a part of the domain. This results in a heterogeneous problem combining an eddy-current model in a subset of the computational domain with a magneto-static model in the remainder of the domain. We adopt a two-domain two-step approach in which the primary variables of the problem are the electric scalar potential and the magnetic vector potential. We show numerical results that validate the formulation.
\end{abstract}
\section{Introduction}

In the simulation of electric arcs in switching devices, electromagnetic fields need to be computed in a domain including moving contacts.
Transient electromagnetic effects are typically neglected in such simulations \cite{Busetto:bib:ArcSimLV}.
However, it has been shown \cite{Busetto:bib:Eddy_current_effects,Busetto:bib:eddy_currents_arc_motion} that eddy currents can have a sizeable effect on the Lorentz forces acting on the plasma arc.
More precisely, the effect of eddy currents in ferromagnetic splitter plates is the largest due to the high magnetic permeability, while in the copper contacts it is much smaller, and in the arc plasma it is negligible due to the low conductivity \cite{Busetto:bib:Eddy_current_effects}.

We conclude that the largest effects of eddy currents are expected in the non-moving metallic parts of the switching device, whereas in the moving part of the domain, eddy currents are negligible.
At the same time, the solution of the equations is more complicated to implement in the moving part which is remeshed, because including the time-dependent eddy current term requires time consuming interpolation.
In this paper, we therefore propose a two-step method coupling the eddy current equations in part of the domain with the magneto-static equations in its complement.

The structure of the paper is as follows.

In Section \ref{Busetto:sec:problem_statement}, we present a coupled potential formulation based on the electric scalar potential $\varphi$ and the magnetic vector potential $\vec{A}$ as primary variables.
Then, we decouple the potentials using two connected ideas inspired by \cite{Busetto:bib:Frequency_stable}: the DC-conduction gauge and the electric circuit element (ECE) boundary conditions.
The DC-conduction gauge allows to recast the original formulation as a two-step algorithm in which we first solve for $\varphi$ (DC-conduction step) and then we use the obtained results to solve for $\vec{A}$ (EC-correction step).
The ECE boundary conditions allow us to correctly impose the boundary conditions in the two-step algorithm.

In Section \ref{Busetto:sec:Weak_formulation} we present a weak formulation suitable for current excitation of the resulting problem.
We also underline an important aspect related to the computation of the voltage that is of interest for engineering applications.

Finally, in Section \ref{Busetto:sec:Numerical_results} we numerically validate our approach through numerical experiments in a conductive cylinder surrounded by air.
First, we consider the cylinder to be made of iron and perform a convergence analysis.
Then, we subdivide the cylinder in three portions in series with the middle one having negligible eddy currents.
We validate the two-domain approach by showing similar results when neglecting eddy currents in the middle portion.

\section{Problem statement}\label{Busetto:sec:problem_statement}
We consider the eddy current equations in the time domain $[0,T] \subset \bbbr$ on a open bounded domain $\Omega \subset \bbbr^3$:
\begin{eqnarray}
  \nabla\times{\vec{E}} + \partial_t\vec{B} & = & \vec{0}, \label{Busetto:eq:max_curl_e} \\
	\nabla\cdot{\vec{B}} & = & 0, \label{Busetto:eq:max_div_b} \\
	\nabla\times{\vec{H}} & =  & \vec{j}, \label{Busetto:eq:max_curl_h}
\end{eqnarray}
and the constitutive relations
\begin{eqnarray}
\vec{B} = \mu \vec{H}, \ \vec{j} = \sigma \vec{E}, \label{Busetto:eq:constitutive_relations}
\end{eqnarray}
where $\vec{E}$ denotes the electric field, $\vec{B}$ the magnetic field, $\vec{H}$ the magnetizing field, and $\vec{j}$ the current density, while the magnetic permeability  $\mu$ and the electrical conductivity $\sigma$ are scalar-valued functions in space and time.

\subsection{Two-domain approach}
In the proposed formulation, we consider eddy currents to be present only in a part of the conductive domain. Therefore, we partition the domain $\Omega$ into the eddy-current subdomain $\Omega_{\mathrm{E}}$ and the magneto-static subdomain $\Omega_{\mathrm{M}}$. The former refers to the part of the domain where eddy fields are present ($\nabla\times{\vec{E}} \neq \vec{0}$) and we also include the insulating subdomain $\Omega_{\mathrm{I}}$ having zero conductivity ($\sigma = 0$), i.e., $\Omega_{\mathrm{I}} \subset \Omega_{\mathrm{E}}$. The latter refers to the part of the conductive domain ($\sigma > 0$) where eddy currents are not taken into account ($\nabla\times{\vec{E}} = \vec{0}$). Furthermore, we define $\Gamma_{\mathrm{I}} := \partial\Omega_{\mathrm{E}}\cap\partial\Omega_{\mathrm{M}}$. 
 
We choose the electric scalar potential $\varphi$ and the magnetic vector potential $\vec{A}$ as primary variables and we proceed as follows.
We set $\vec{B} = \nabla\times{\vec{A}}$, so that (\ref{Busetto:eq:max_div_b}) is satisfied by construction.
In $\Omega_{\mathrm{M}}$ we neglect $\partial_t\vec{B}$ in (\ref{Busetto:eq:max_curl_e}), so that $\nabla\times{\vec{E}} = \vec{0}$. Therefore, we can set $\vec{E} = -\nabla{\varphi}$ in $\Omega_{\mathrm{M}}$, assuming $\Omega_{\mathrm{M}}$ to be simply connected. Moreover, we set $\vec{E} = - \partial_t\vec{A} - \nabla{\varphi}$ in $\Omega_{\mathrm{E}}$, assuming $\Omega_E$ to be simply connected.
Therefore, (\ref{Busetto:eq:max_curl_h}) is the only remaining equation that we need to solve, all others being satisfied by construction.

Given the previous considerations and the constitutive relations (\ref{Busetto:eq:constitutive_relations}), we can write (\ref{Busetto:eq:max_curl_h}) as a function of the primary variables $\varphi$ and $\vec{A}$ as follows:
\begin{eqnarray}
  \sigma\partial_t\vec{A} + \sigma\nabla{\varphi} + \nabla\times\left(\mu^{-1}\nabla\times{\vec{A}}\right) & = & \vec{0} \quad\mbox{in }\Omega_{\mathrm{E}}, \label{Busetto:eq:ec_curl}\\
  \sigma\nabla{\varphi} + \nabla\times\left(\mu^{-1}\nabla\times{\vec{A}}\right) & = & \vec{0} \quad\mbox{in }\Omega_{\mathrm{M}}.\label{Busetto:eq:ms_curl} 
  \end{eqnarray}
We choose interface conditions to act on the two potentials separately and to ensure continuity of the normal current density and normal magnetic flux density:
\begin{eqnarray}
	\left[\varphi\right]_{\Gamma_{\mathrm{I}}} = \vec{0}, \quad \left[\vec{n}\times\vec{A}\right]_{\Gamma_{\mathrm{I}}} = \vec{0}, \quad\left[\vec{n}\times\vec{H}\right]_{\Gamma_I} & = & \vec{0}\quad\mbox{on }\Gamma_{\mathrm{I}},\label{Busetto:eq:interface_cond}
\end{eqnarray}
where $\left[\psi\right]_{\Gamma_I} = \psi\vert_{\Gamma_{E}}\vec{n}_E + \psi\vert_{\Gamma_M}\vec{n}_M$ represents the jump of a scalar function $\psi$ across $\Gamma_I$ where $E$ and $M$ stand for the two sides of an interface $\Gamma_I$ and $\vec{n}_k$ denotes the outer normal to $\Omega_k$ on $\Gamma_I$.
We denote by $\left[\vec{n}\times\vec{\psi}\right]=\vec{n}_E\times\vec{\psi}\vert_{\Gamma_E}+\vec{n}_M\times\vec{\psi}\vert_{\Gamma_M}$ the (rotated) jump of the tangential component of a vector-valued function $\vec{\psi}$.
We notice that this choice implies that $\left[\vec{n}\times\vec{E}\right]_{\Gamma_{\mathrm{I}}}=-\partial_t(\vec{n}_E\times\vec{A})\neq \vec{0}$, i.e., in general the tangential component of $\vec{E}$ is not continuous at the intersection $\Gamma_{\mathrm{I}}$. We point out that in the insulating subdomain $\Omega_{\mathrm{I}}$ (\ref{Busetto:eq:ec_curl}) simplifies to  $\nabla\times\left(\mu^{-1}\nabla\times{\vec{A}}\right)  =  \vec{0}$. Moreover, we remark that in the special case $\Omega_{\mathrm{M}} = \emptyset$ the problem reduces to the classical eddy current formulation.

\subsection{Two-step approach}
In what follows we decouple the potential exploiting the DC-conduction gauge \cite{Busetto:bib:Frequency_stable} and the electric circuit element (ECE) boundary conditions \cite{Busetto:bib:Electric_Circuit_Element}.
The possibility to decouple the potentials relies on the freedom to choose $\varphi$ inside $\Omega_{\mathrm{E}}$ (gauge freedom).
Therefore, we choose $\varphi$ as the solution of the equation obtained by taking the divergence of (\ref{Busetto:eq:ms_curl}) and extending it to the whole domain $\Omega$, i.e.,
\begin{eqnarray}
  \nabla\cdot{\left(-\sigma\nabla{\varphi}\right)} & = & 0 \qquad\mbox{in }\Omega. \label{Busetto:eq:dc_cond}
\end{eqnarray}
Equation (\ref{Busetto:eq:dc_cond}) is known as the DC-conduction equation. 

Given this, we proceed applying a sequential two-step approach.
In the first step (DC-conduction) we solve (\ref{Busetto:eq:dc_cond}) and we compute $\varphi$.
Then we exploit this result as a gauge condition for the second step (EC-correction) in which we solve the following equations:
\begin{eqnarray}
  \sigma\partial_t\vec{A} + \nabla\times\left(\mu^{-1}\nabla\times{\vec{A}}\right) & =  & -\sigma\nabla{\varphi} \qquad\mbox{in }\Omega_{\mathrm{E}}, \label{Busetto:eq:ec_corr} \\
	\nabla\times\left(\mu^{-1}\nabla\times{\vec{A}}\right) & =  & -\sigma\nabla{\varphi} \qquad\mbox{in }\Omega_{\mathrm{M}},\label{Busetto:eq:ms_corr}
\end{eqnarray}
complemented by the interface conditions (\ref{Busetto:eq:interface_cond}).
The second step has the two-fold purpose to correct the electric field introducing the inductive effects and to compute the magnetic field. 

\subsection{ECE boundary conditions}\label{Busetto:sec:ECE_boundary_conditions}
We need to impose suitable boundary conditions that enable coupling with circuit models and that do not couple the electric and magnetic potentials back again.
We choose the electric-circuit-element (ECE) boundary conditions \cite{Busetto:bib:Electric_Circuit_Element}.
We assume that the boundary $\partial \Omega$ with outer unit normal $\vec{n}$ is partitioned into two non-touching contacts or electric ports, $\Pi_1$ and $\Pi_2$ and the insulated boundary part $\Pi_0$. The ECE boundary conditions approach is based on the following assumptions:
\begin{enumerate}
\item \label{Busetto:item:ece_ass_1}no inductive coupling with the exterior: $(\partial_t\vec{B}) \cdot \vec{n} =  0 \Leftrightarrow (\nabla\times{\vec{E}}) \cdot \vec{n} = 0$ on $\partial \Omega$;
\item \label{Busetto:item:ece_ass_2}electric ports are equipotential surfaces: $\left(\vec{n}\times\vec{E}\right)\times\vec{n} = \vec{0}$ on $\Pi_j$, $j = 1,2$;
\item \label{Busetto:item:ece_ass_3}no electric currents can enter the insulated boundary part: $(\nabla\times \vec{H}) \cdot \vec{n} = 0 \Leftrightarrow \vec{E} \cdot \vec{n} = 0$ on $\Gamma_0$.
\end{enumerate}

Notice that from assumptions \ref{Busetto:item:ece_ass_1} and \ref{Busetto:item:ece_ass_2} we conclude that
\begin{eqnarray*}
  \left(\vec{n}\times\vec{E}\right)\times\vec{n} & = & \left(\vec{n}\times\nabla{\varphi_{\mathrm{tot}}}\right)\times\vec{n}\quad\mbox{on }\partial\Omega,\\
\mbox{for some}\quad\varphi_{\mathrm{tot}} \in H^1_{\mathrm{ece}}(\Omega) & := & \left\{\phi\in H^1(\Omega) : \phi\vert_{\Pi_j}=\mbox{const}, \ j = 1,2\right\}.
\end{eqnarray*}

The ECE boundary conditions are applied in both the DC-conduction and the EC-correction steps.
In the DC-conduction step, $\varphi$ computed in (\ref{Busetto:eq:dc_cond}) satisfies assumption 1. Moreover, by choosing $\varphi: [0,T] \rightarrow H^1_{\mathrm{ece}}(\Omega)$, also assumption 2 is automatically satisfied. In order to satisfy also assumption 3, we complement by the zero-current condition $-\sigma\nabla{\varphi}\cdot\vec{n}=0$ on the insulating part $\Gamma_0$.

In the EC-correction step (\ref{Busetto:eq:ec_corr}), in order to satisfy assumptions 1 and 2, we choose $\vec{A}: [0,T] \rightarrow H_{\mathrm{ece}}(\mbox{curl},\Omega)$ where
\begin{eqnarray*}
  H_{\mathrm{ece}}(\mbox{curl},\Omega) := \left\{ \vec{a} \in H(\mbox{curl},\Omega) : \vec{n}\times\vec{a}=\vec{n}\times\nabla{\eta}\mbox{ on }\partial\Omega,\ \eta\in H^1_{\mathrm{ece}}(\Omega) \right\}. 
\end{eqnarray*}
Indeed, this is equivalent to impose $\nabla\times \vec{A} \cdot \vec{n} = 0$ on $\partial \Omega$ and $\left(\vec{n}\times\vec{A}\right)\times\vec{n} = \vec{0}$ on $\Pi_j$,  $j = 1,2$. Finally, in order to satisfy assumption 3, we complement by the zero-current condition $\nabla\cdot (\frac{1}{\mu} \nabla\times \vec{A} \times \vec{n} )$ on $\Gamma_0$.

{Given that $\vec{E} = - \partial_t\vec{A} - \nabla{\varphi}$ in $\Omega_{\mathrm{E}}$, in case eddy currents are present in the whole conductive domain, i.e., $\Omega_{\mathrm{M}} = \emptyset$, the total potential $\varphi_{\mathrm{tot}}$ at a port $\Pi_j$ is given by the superposition of the potentials $V_{j,\mathrm{DC}}$ and $V_{j,\mathrm{EC}}$ obtained in the two steps, i.e.,
\begin{eqnarray}\label{Busetto:eq:total_voltage}
  V_j := \varphi_{\mathrm{tot}}\vert_{\Pi_j}  =  \varphi\vert_{\Pi_j} + \partial_t\eta\vert_{\Pi_j} =  V_{j,\mathrm{DC}} + {V_{j,\mathrm{EC}}}, \ \ \ \eta: [0,T] \rightarrow H^1_{\mathrm{ece}}(\Omega).
\end{eqnarray}}
Similarly, the total current $I_j$ flowing through port $\Pi_j$ is given by the superposition of the DC-current $I_{\mathrm{DC}}$ and the EC-current $I_{\mathrm{EC}}$ both accessible through flux integrals:
\begin{eqnarray}
  I_j & = & \int_{\Pi_j}-\sigma\nabla{\varphi}\cdot\vec{n} \ \mbox{d} S + \int_{\Pi_j}-\sigma\partial_t\vec{A}\cdot\vec{n} \ \mbox{d} S \nonumber\\
      & = & \int_{\Pi_j}\vec{j}_{\mathrm{DC}}\cdot\vec{n} \ \mbox{d} S + \int_{\Pi_j}\vec{j}_{\mathrm{EC}}\cdot\vec{n} \ \mbox{d} S
				=   I_{j,\mathrm{DC}} + I_{j,\mathrm{EC}}.\label{Busetto:eq:total_current}
\end{eqnarray}

We impose the full currents and voltages in the DC-conduction step as in \cite{Busetto:bib:Frequency_stable}.
Indeed, in $\Omega_{\mathrm{M}}$, only $\varphi$ can represent the total current, so if we were to impose zero currents and voltages in the first step, we would have $\varphi=0 \ \forall (\vec{x},t)$ and never get any current in $\Omega_{\mathrm{M}}$.
Hence the full boundary conditions are applied in the first step, whereas the second step is computed with the same type of boundary conditions of the full step, but with their values set to zero.
We prescribe a current $I_1$ at port $\Pi_1$ and fix the potential $V_2$ at $\Pi_2$ for uniqueness of $\varphi$. In the DC-conduction step fix $I_{1,\mathrm{DC}} := I_1$ and $V_{2,\mathrm{DC}} := V_2$ and in the EC-correction step set $I_{1,\mathrm{EC}} := 0$ and $V_{2,\mathrm{EC}} := 0$.
For simplicity of presentation and without loss of generality, we set $V_{2,\mathrm{DC}}=V_2=0$.

\section{Weak formulation and discretization}\label{Busetto:sec:Weak_formulation}
In the following we state the weak formulation for both steps of the considered problem.

\subsection{DC-conduction step}
The ECE boundary conditions are captured representing $\varphi$ as
\begin{eqnarray*}
  \varphi = \tilde{\varphi} + V_{1,\mathrm{DC}} \Phi_1,
\end{eqnarray*}
where $\tilde{\varphi}: [0,T] \rightarrow H^1_{\mathrm{ece},0}(\Omega)$, $H^1_{\mathrm{ece},0}(\Omega) := \left\{\phi\in H^1(\Omega) : \phi\vert_{\Pi_j}= 0, \ j = 1,2\right\}$, $\Phi_i \in H^1_{\mathrm{ece}}(\Omega)$ with $\Phi_i\vert_{\Pi_j} = \delta_{ij}$, $i,j\in\{1,2\}$ and $V_{1,DC}: [0,T] \rightarrow \bbbr$

Then the variational formulation reads as follows.

Set $I_{1,\mathrm{DC}} = I_1$ and seek $\tilde{\varphi}: [0,T] \rightarrow H^{1}_{\mathrm{ece},0}(\Omega)$, $V_{1,\mathrm{DC}}: [0,T] \rightarrow \bbbr$ such that $\forall \tilde{\varphi}' \in H^{1}_{\mathrm{ece},0}(\Omega)$, $V_1'\in\bbbr$
\begin{eqnarray}\label{Busetto:eq:dc_cond_weak}
(\sigma \nabla{{(\tilde\varphi} + V_{1,\mathrm{DC}} \Phi_1)}, \nabla{ \tilde{\varphi}'} + V_1'\nabla{\Phi_1})_{\Omega}  & = & -I_{1,\mathrm{DC}} V_1'.
\end{eqnarray}

\subsection{EC-correction step}
The ECE boundary conditions are captured representing $\vec{A}$ as
\begin{eqnarray*}
\vec{A} = \vec{\tilde{A}} + W_{1,\mathrm{EC}} \nabla{\Phi_1} + W_{2,\mathrm{EC}}\nabla\Phi_2,
\end{eqnarray*}
where $W_{j,\mathrm{EC}}: [0,T] \rightarrow \bbbr$, $\partial_t W_{j,\mathrm{EC}}=V_{j,\mathrm{EC}}$, $j\in\{1,2\}$, and $\vec{\tilde{A}}: [0,T] \rightarrow H_{\mathrm{ece},0}{(\mbox{curl},\Omega)}$ with
\begin{eqnarray*}
H_{\mathrm{ece},0}(\mbox{curl},\Omega) = \{ \vec{v} \in H_{\mathrm{ece}}(\mbox{curl},\Omega), \exists \eta \in H^1_{\mathrm{ece},0}(\Omega): \vec{n} \times \vec{v}\vert_{\partial \Omega}  =  \vec{n} \times \nabla{\eta}\vert_{\partial \Omega} \}. 
\end{eqnarray*}

We multiply (\ref{Busetto:eq:ec_corr}) with a test function  ${\vec{A}'} =  \vec{\tilde{A}}'+ W_{1}' \nabla{\Phi_1} + W_{2}' \nabla{\Phi_2}$, with $\vec{\tilde{A}}' \in H_{ece,0}(\mbox{curl},\Omega)$, $W_1', W_2' \in \bbbr$ and we apply the Green's formula on $\Omega_{\mathrm{E}}$
\begin{eqnarray*}
        \left(\sigma\partial_t\vec{A},\vec{A}'\right)_{\Omega_{\mathrm{E}}} + 
	      \left(\mu^{-1}\nabla\times{\vec{A}},\nabla\times{\vec{A}'}\right)_{\Omega_{\mathrm{E}}} +
        \left\langle \vec{n}\times \left(\mu^{-1}\nabla\times{\vec{A}}\right) , \vec{A}'\right\rangle_{\partial\Omega_{\mathrm{E}}}&&\\
  = \left(-\sigma\nabla{\varphi},\vec{A}'\right)_{\Omega_{\mathrm{E}}}.&&
\end{eqnarray*}

Then we multiply (\ref{Busetto:eq:ms_corr}) with the same test function $\vec{{A}'}$ and we apply the Green's formula on $\Omega_{\mathrm{M}}$
\begin{eqnarray*}
 \left(\mu^{-1}\nabla\times{\vec{A}},\nabla\times{\vec{A}'}\right)_{\Omega_{\mathrm{M}}} 
 + \left\langle \vec{n}\times \left(\mu^{-1}\nabla\times{\vec{A}}\right) , \vec{A}' \right\rangle_{\partial\Omega_{\mathrm{M}}}
     & = & \left(-\sigma\nabla{\varphi},\vec{A}'\right)_{\Omega_{\mathrm{M}}}.
\end{eqnarray*}
Adding the previous two equations together and noting that unit normal $\vec{n}$ on $\Gamma_I$ has opposite direction in $\partial \Omega_M$ and in $\partial \Omega_E$,  we obtain
\begin{eqnarray*}
        &&\left(\sigma\partial_t\vec{A},\vec{A}'\right)_{\Omega_{\mathrm{E}}} +
	      \left(\mu^{-1}\nabla\times{\vec{A}},\nabla\times{\vec{A}'}\right)_{\Omega} + \\
	     && \left\langle \left[\vec{n}\times \left(\mu^{-1}\nabla\times{\vec{A}}\right)\right],\vec{A}'\right\rangle_{\Gamma_{\mathrm{I}}} +
				\left\langle \vec{n}\times \left(\mu^{-1}\nabla\times{\vec{A}}\right),\vec{A}'\right\rangle_{\partial\Omega}
  =  \left(-\sigma\nabla{\varphi},\vec{A}'\right)_{\Omega}.
\end{eqnarray*}

As $\left[\vec{n}\times\vec{H}\right]_{\Gamma_I}=\vec{0}$ according to 
(\ref{Busetto:eq:interface_cond}), the jump term on $\Gamma_{\mathrm{I}}$ drops.
 We notice that $\left(\mu^{-1}\nabla \times {\vec{A}},\nabla \times {{\vec{{A}}'}}\right)_{\Omega}$ simplifies to $\left(\mu^{-1}\nabla \times {\vec{\tilde{A}}},\nabla \times {{\vec{\tilde{A}}'}}\right)_{\Omega}$.
Moreover, noting that $\mu^{-1} \nabla \times {\vec{A}} = \mu^{-1} \vec{B} = \vec{H}$, that $\vec{\tilde{A}}'\times\vec{n}=\vec{0}$ on the boundary, integrating by parts, and using $\nabla \cdot{\vec{n}\times{\vec{H}}}=\vec{H}\cdot\nabla\times{\vec{n}}-\vec{n}\cdot\nabla\times{\vec{H}}$, $\nabla\times{\vec{n}}=\vec{0}$, $\nabla\times{\vec{H}} = \vec{j}$ and the zero boundary condition of Section \ref{Busetto:sec:ECE_boundary_conditions}, the boundary term can be treated as follows
\begin{eqnarray*}
      \left\langle \vec{n}\times \vec{H}, \vec{A}^{'}\right\rangle_{\partial \Omega}
    = I_1 W_{1}' +  I_2 W_{2}'.
\end{eqnarray*}
Furthermore, noting that $I_j-I_{j,DC} = I_{j,EC}$, we have
\begin{eqnarray*}
  \left(-\sigma\nabla{\varphi},{\vec{{A}}'}\right)_{\Omega} - \left( I_1 W_{1}' +  I_2 W_{2}'\right) 
  & = & \left(-\sigma\nabla{\varphi},\vec{\tilde{A}}'\right)_{\Omega} - \left( I_{1,EC} W_{1}' +  I_{2,EC} W_{2}'\right).
\end{eqnarray*}
The variational formulation reads as follows:
For $W_{2,\mathrm{EC}} = 0\Rightarrow \partial_t W_{2,\mathrm{EC}} = V_{2,\mathrm{EC}} = 0$ and $I_{1,\mathrm{EC}} = 0$, seek $\vec{\tilde{A}}: [0,T] \rightarrow H_{\mathrm{ece},0}{(\mbox{curl},\Omega)}$, $W_{1,\mathrm{EC}}: [0,T] \rightarrow \bbbr$ such that $\forall \vec{\tilde{A}}' \in H_{\mathrm{ece},0}(\mbox{curl},\Omega)$, $W_1'\in\bbbr$
\begin{eqnarray}\label{Busetto:eq:a_weak_current}
        \left(\sigma\partial_t(\vec{\tilde{A}} + W_{1,\mathrm{EC}} \nabla \Phi_1),\vec{\tilde{A}}' + W_1'\nabla \Phi_1\right)_{\Omega_{\mathrm{E}}}+ \left(\mu^{-1}\nabla\times{\vec{\tilde{A}}},\nabla\times{\vec{\tilde{A}}'}\right)_{\Omega} &&\nonumber \\
   = \left(-\sigma\nabla{\varphi},\vec{\tilde{A}}'\right)_{\Omega}.&&
\end{eqnarray}

\subsection{Voltage reconstruction} \label{Busetto:sec:VoltageReconstruction}

In case eddy currents are present in the whole conductive domain, i.e. $\Omega_{\mathrm{M}} = \emptyset$, the voltage drop can be computed as $U_{1,2} = V_2 - V_1$, with $V_1$ being the total potential given by $V_1=V_{1,\mathrm{DC}}+V_{1,\mathrm{EC}}$, where
\begin{equation}\label{Busetto:eq:voltage_ec}
    V_{1,\mathrm{EC}} = \partial_t W_{1,\mathrm{EC}}.
\end{equation}

In this case, it has been proved in \cite{Busetto:bib:Frequency_stable} that, the total power $P_{\mathrm{total}}$ transferred to the system and computed as the sum of ohmic power $P_{\mathrm{Ohm}}$ and magnetic power $P_{\mathrm{mag}}$ equals the total power computed as the product of the total current $I_{1}$ given by (\ref{Busetto:eq:total_current}) and the voltage drop $U_{1,2}$ computed from (\ref{Busetto:eq:total_voltage}), i.~e.,
\begin{eqnarray*}\label{Busetto:eq:total_power}
 P_{\mathrm{total}}  = P_{\mathrm{Ohm}} + P_{\mathrm{mag}} = U_{1,2} I_1,
\end{eqnarray*}
where 
\begin{eqnarray*}
P_{\mathrm{Ohm}} & = & (\vec{j},\vec{E})_{\Omega} = (\sigma\nabla \varphi, \nabla \varphi)_{\Omega_{\mathrm{M}}} + (\sigma(\nabla\varphi+\partial_t\vec{A}),\nabla\varphi+\partial_t\vec{A})_{\Omega_{\mathrm{E}}}, \label{Busetto:eq:POhm}\\
P_{\mathrm{mag}} & = & (\vec{H},\partial_t\vec{B})_{\Omega} = ({\mu}^{-1} \nabla\times{\vec{A}},\nabla\times{\partial_t\vec{A}})_{\Omega}. \label{Busetto:eq:Pmag}
\end{eqnarray*}

In case eddy currents are not present in the whole conductive domain, i.e., $\Omega_{\mathrm{M}} \neq \emptyset$, it can be shown that the quantity $U_{1,2} := V_2-V_1$ computed as above will not coincide with the actual voltage drop (see the numerical results of Section \ref{Busetto:sec:Test_cases_2_and_3}).
Hence $P_{\mathrm{total}} = P_{\mathrm{Ohm}} + P_{\mathrm{mag}} \neq U_{1,2}I_1$, because in this case (\ref{Busetto:eq:total_voltage}) does not allow to properly reconstruct the voltages $V_j$.
Indeed, $V_{j,\mathrm{EC}}$ is no longer uniquely defined and therefore does not allow to compute the inductive correction to the voltage.
However, the actual voltage drop $U_{1,2}$ can be reconstructed a posteriori from the computed total power $P_{\mathrm{total}}$ as $U_{1,2}=P_{\mathrm{total}}/I_1$.
As there is only one power, we can only compute one voltage, which is what restricts us to two ports.

\subsection{Discretization}\label{Busetto:sec:Discretization}

We discretize (\ref{Busetto:eq:a_weak_current}) and (\ref{Busetto:eq:voltage_ec}) in time. We initialize with zero fields, perform one step of the backward Euler scheme and then we apply the BDF2 scheme.

For the space discretization of (\ref{Busetto:eq:dc_cond_weak}) and (\ref{Busetto:eq:a_weak_current}), we adopt the standard first order Lagrangian nodal finite element space for the approximation of $\varphi$, $\Phi_1$ and $\Phi_2$, and the standard first order N\'ed\'elec edge finite element space for $\vec{\tilde{A}}$.

The resulting symmetric positive semi-definite system of linear equations for $\vec{\tilde{A}}$ is solved applying the auxiliary space preconditioner for edge elements \cite{Busetto:bib:Nodal_auxiliary}.

\section{Numerical results}\label{Busetto:sec:Numerical_results}
In this section we perform numerical experiments to assess the validity of the described theoretical approach. We consider a domain $\Omega$ consisting of an inner conductive cylinder ($\sigma > 0$) surrounded by a non-conductive cylinder consisting of air ($\sigma = 0$). The air subdomain coincides with $\Omega_{\mathrm{I}}$, whereas the inner conductive cylinder coincides with $(\Omega_{\mathrm{E}} \backslash \Omega_{\mathrm{I}}) \cup \Omega_{\mathrm{M}}$.
We focus on three different problems:
\begin{itemize}
\item Test case 1 - homogeneous inner conductive cylinder made of iron with eddy currents present everywhere ($\Omega_{\mathrm{M}} = \emptyset$);
\item Test cases 2 and 3 - inner conductive cylinder split in three different portions in series consisting of iron (first and third portions), copper (second portion):\\
- Test case 2 - eddy currents present everywhere ($\Omega_{\mathrm{M}} = \emptyset$)\\
    - Test case 3 - eddy currents present only in the iron part ($\Omega_{\mathrm{M}} \neq \emptyset$).
\end{itemize}
The imposed current excitation is sinusoidal in time, i.e., $I_1(t) = I_0 \cos (\omega t + \pi)$.

\subsection{Choice of physical parameters}\label{Busetto:sec:ChoiceOfPhysicalParameters}

We choose a frequency $f=50\mbox{ Hz}$ ($\omega=2\pi f$), and the following values for the material parameters: for iron an electric conductivity $\sigma_{\mathrm{Fe}} = 10^7\mbox{ S/m}$ and a relative magnetic permeability $\mu_{\mathrm{r,Fe}} = 1500$; for copper, $\sigma_{\mathrm{Cu}} = 6\cdot 10^7\mbox{ S/m}$ and $\mu_{\mathrm{r,Cu}} = 1$.

We choose the radius $R$ of the inner cylinder between the skin depths of iron and copper, such that there is a pronounced skin effect in iron and near uniform current density in copper.
The skin depth $\delta$ of a sinusoidal (in time) current along an infinite half-plane is $\delta = \sqrt{2/(\mu\sigma\omega)}$, which yields $\delta_{\mathrm{Fe}} = 0.41\mbox{ mm}$ for iron and $\delta_{\mathrm{Cu}} = 6.50\mbox{ mm}$ for copper.
Therefore we choose a radius $R=3\mbox{ mm}$ for the inner cylinder and an external radius of the domain $R_{\mathrm{air}}=8\mbox{ mm}$ for all test cases.

\subsection{Test case 1}
This benchmark problem has an analytic solution which allows us to perform a convergence analysis.
We remark that the analytic solution holds for an infinite cylinder.
However, the exact solution on a cylinder of finite length equipped with boundary conditions satisfied by the analytic solution coincides with the analytic solution.

\subsubsection{Analytic Solution}\label{Busetto:sec:Analytic_solution}
For a given total current, expressed as a complex valued phasor $I$, the analytic solution of the (axial) current density $j(r)$ (again as a complex phasor) in (a finite section of) an infinitely long cylinder of radius $R$ depends only on the radial coordinate $r$ and can be computed as \cite{Busetto:bib:Weeks_Transmission}
\begin{eqnarray*}
  j(r) = \frac{k I}{2\pi R}\frac{J_0(kr)}{J_1(kR)}\quad r\leq R; \qquad &&
  j(r) = 0 \quad r>R,
\end{eqnarray*}
where $k=(1-i)/\delta$ is the (complex) wave number in the conductor, and $J_{\alpha}$ is the (complex valued) Bessel function of the first kind of order $\alpha$.
From this can be derived the (circumferential) magnetic field $B(r)$ as 
\begin{eqnarray*}
  B(r) = \frac{I \mu}{2\pi R} \frac{J_1(kr)}{J_1(kR)}\quad r\leq R; \qquad &&
  B(r) = \frac{I\mu}{2 \pi r}\quad r>R.
\end{eqnarray*}

\subsubsection{Numerical results}
We consider four consecutive mesh refinements $\mathcal{T}_j$, $j\in\{1,2,3,4\}$ of the cylindrical domain $\Omega$.
Furthermore, we consider a time interval $\mathcal{I}_{T} = [0,T]$ with $T = 0.14\mbox{ s}$ corresponding to seven periods at $f = 50\mbox{Hz}$.
Indeed, since we initialize with zero fields as explained in Section \ref{Busetto:sec:Discretization}, we need to simulate a sufficient number of periods before the solution approaches a stationary regime.

\begin{figure}[htbp]
\centering
\subfloat[][$L^2$-norm]{
\includegraphics[width=0.48\textwidth]{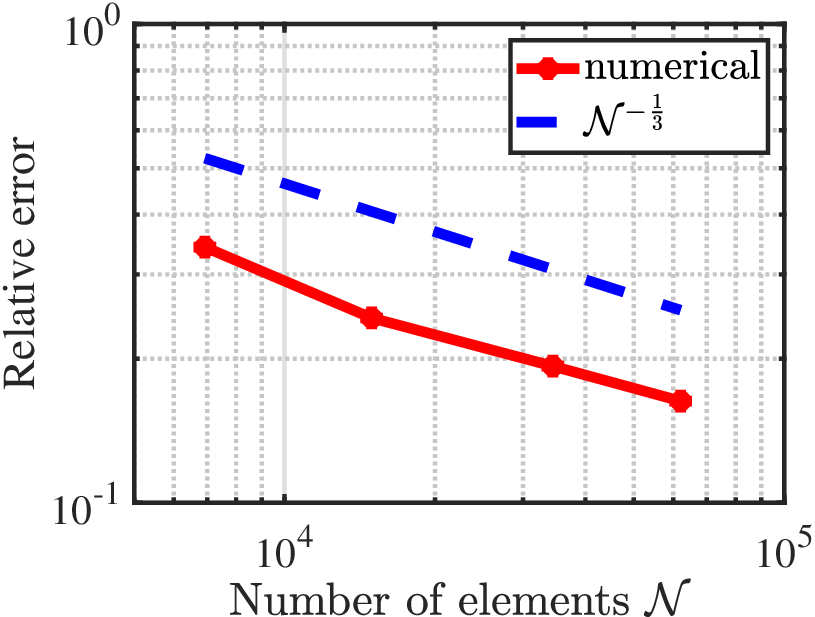}
\label{Busetto:fig:Busetto_current_density_iron_current}}
\subfloat[][$H(\mbox{curl})$-seminorm]{
\includegraphics[width=0.48\textwidth]{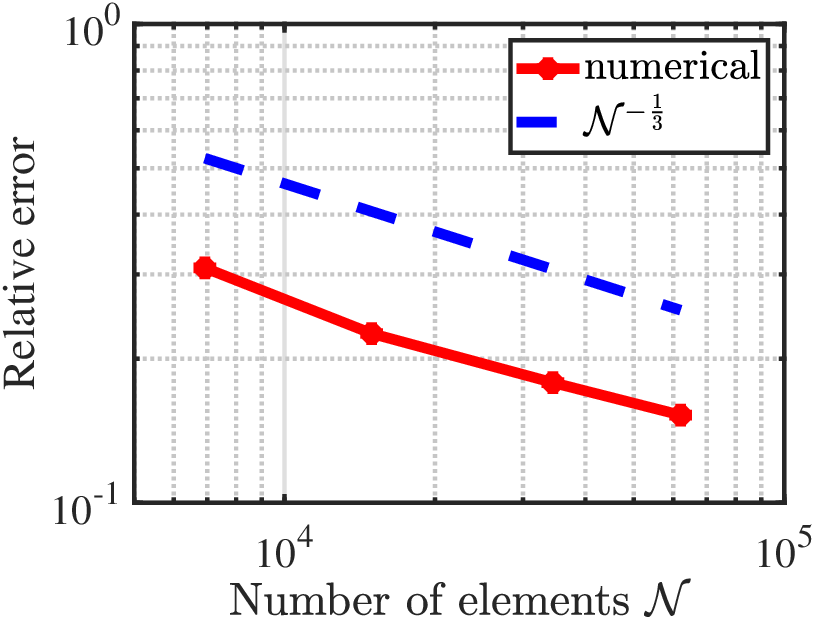}
\label{Busetto:fig:Busetto_magnetic_field_iron_current}}
\caption{ Test case 1: (a) convergence in $L^2$-norm and (b) convergence in $H(\mbox{curl})$-seminorm.}
\label{Busetto:fig:errors_iron_currents}
\end{figure}

From theory on similar problems \cite{Busetto:bib:Monk}, we expect convergence order $O(\mathcal{N}^{-\frac{1}{3}})$ both in $L^2$-norm and in $H(\mbox{curl})$-seminorm, where $\mathcal{N}$ denotes the number of elements.
In Fig. \ref{Busetto:fig:errors_iron_currents} we plot the relative error as the difference between the numerical and the analytic solution in $L^2$-norm and in $H(\mbox{curl})$-seminorm integrated over the last time period for an increasing number of elements.
The slope of the line connecting the numerical values is comparable to the slope of the line representing the expected order of convergence.
Therefore, we conclude that the numerical solutions converge to the analytic solutions with the correct order.

\subsection{Test cases 2 and 3}\label{Busetto:sec:Test_cases_2_and_3}
\begin{figure}[htbp]
\centering
\includegraphics[width=0.8\textwidth]{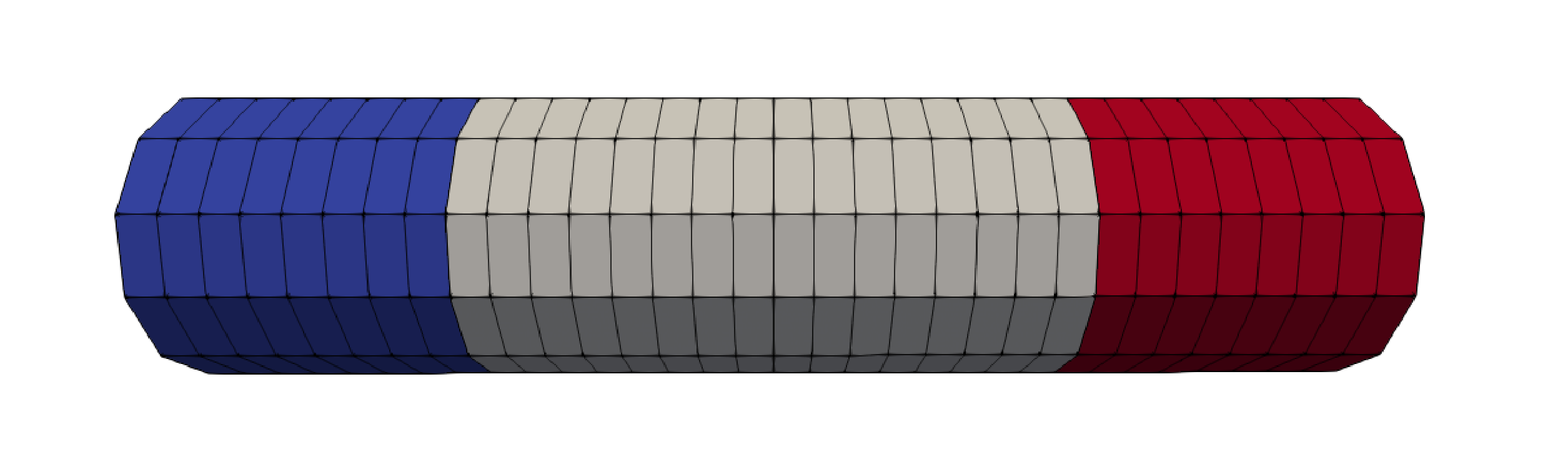}
\caption{Test case 2 and 3 (prisms): Three portion inner cylinder. The blue and the red portions refers to the iron part, whereas the gray portion refers to the copper part.}
\label{Busetto:fig:three_portion_cylinder}
\end{figure}

A representation of the three-portion inner cylinder is shown in Fig. \ref{Busetto:fig:three_portion_cylinder}.
The physical data of the problem are reported in Section \ref{Busetto:sec:ChoiceOfPhysicalParameters}.
Concerning the geometry, we consider five cylinders with increasing lengths.
From one mesh to the subsequent one we double the length of each one of the three portions in series. Moreover, the length of the copper-middle portion is always double with respect to the length of the iron portions. In our analysis, we consider meshes consisting of prisms and layers of fixed size. We consider three consecutive mesh refinements for each one of the given cylinders. 

\begin{figure}[htbp]
\centering
\subfloat[][Iron]{
\includegraphics[width=0.48\textwidth]{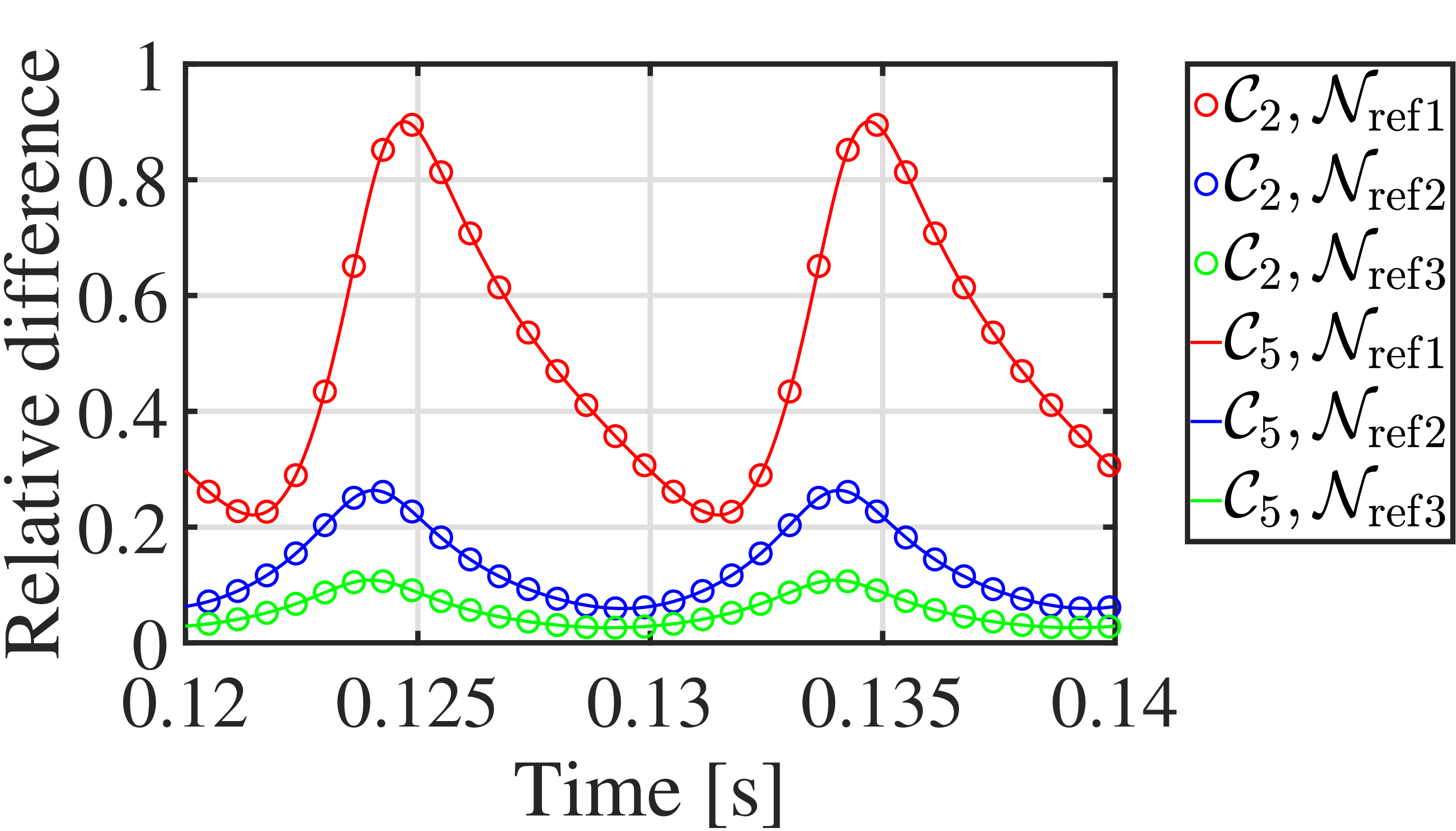}
\label{Busetto:fig:BusettoCurrentDensityIntersection1}}
\subfloat[][Copper]{
\includegraphics[width=0.48\textwidth]{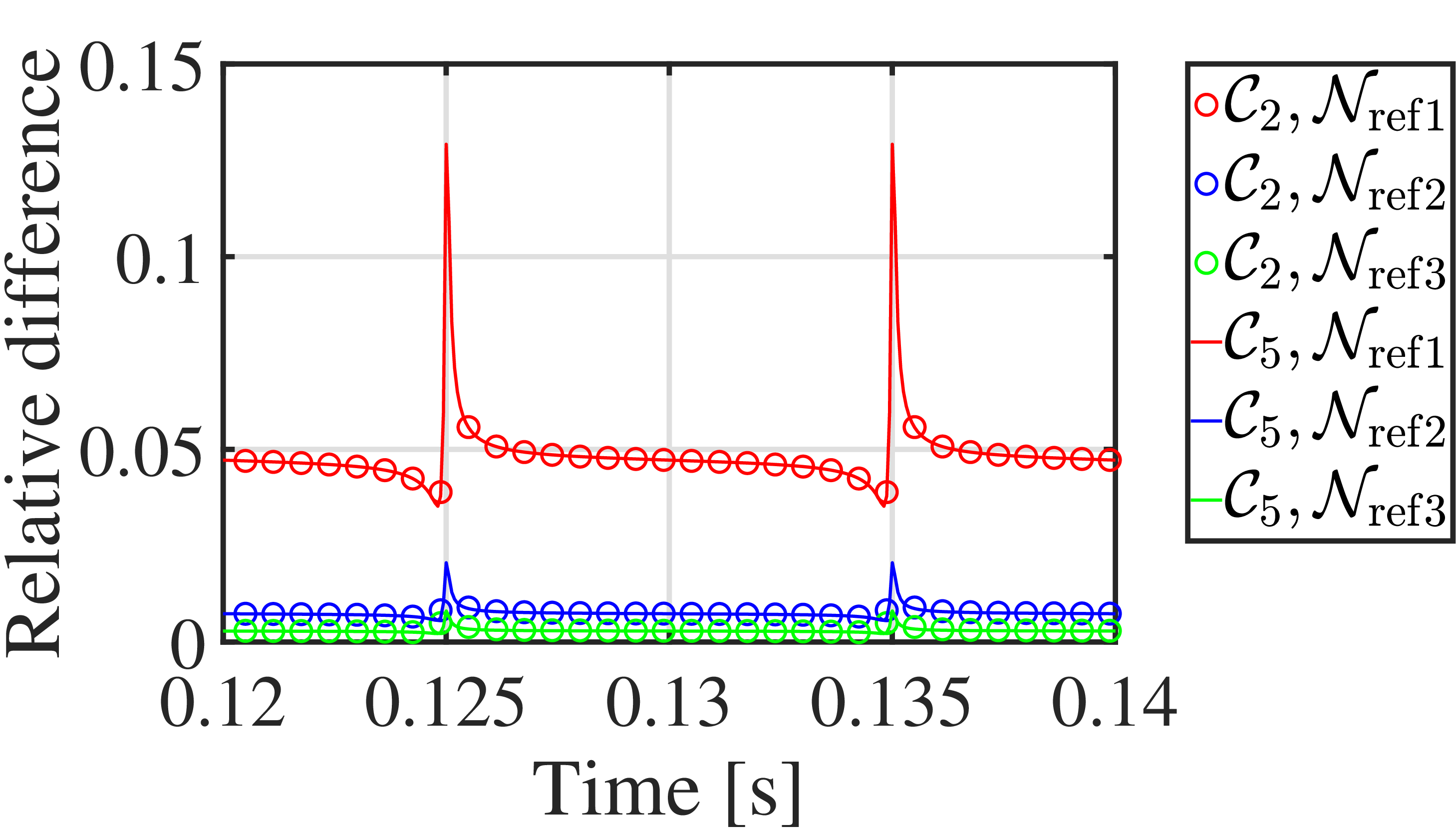}
\label{Busetto:fig:BusettoCurrentDensityIntersection2}}
\caption{Test case 2: Evolution in time of $L^2$-norm of difference in current density between numeric solution on three-portion cylinder and analytic solution on infinite cylinder: (a) in first (or third) portion of inner cylinder (iron) and (b) in second portion of inner cylinder (copper).}
\label{Busetto:fig:BusettoCurrentDensityIntersections}
\end{figure}

\begin{figure}[htbp]
\centering
\subfloat[][Iron]{
\includegraphics[width=0.48\textwidth]{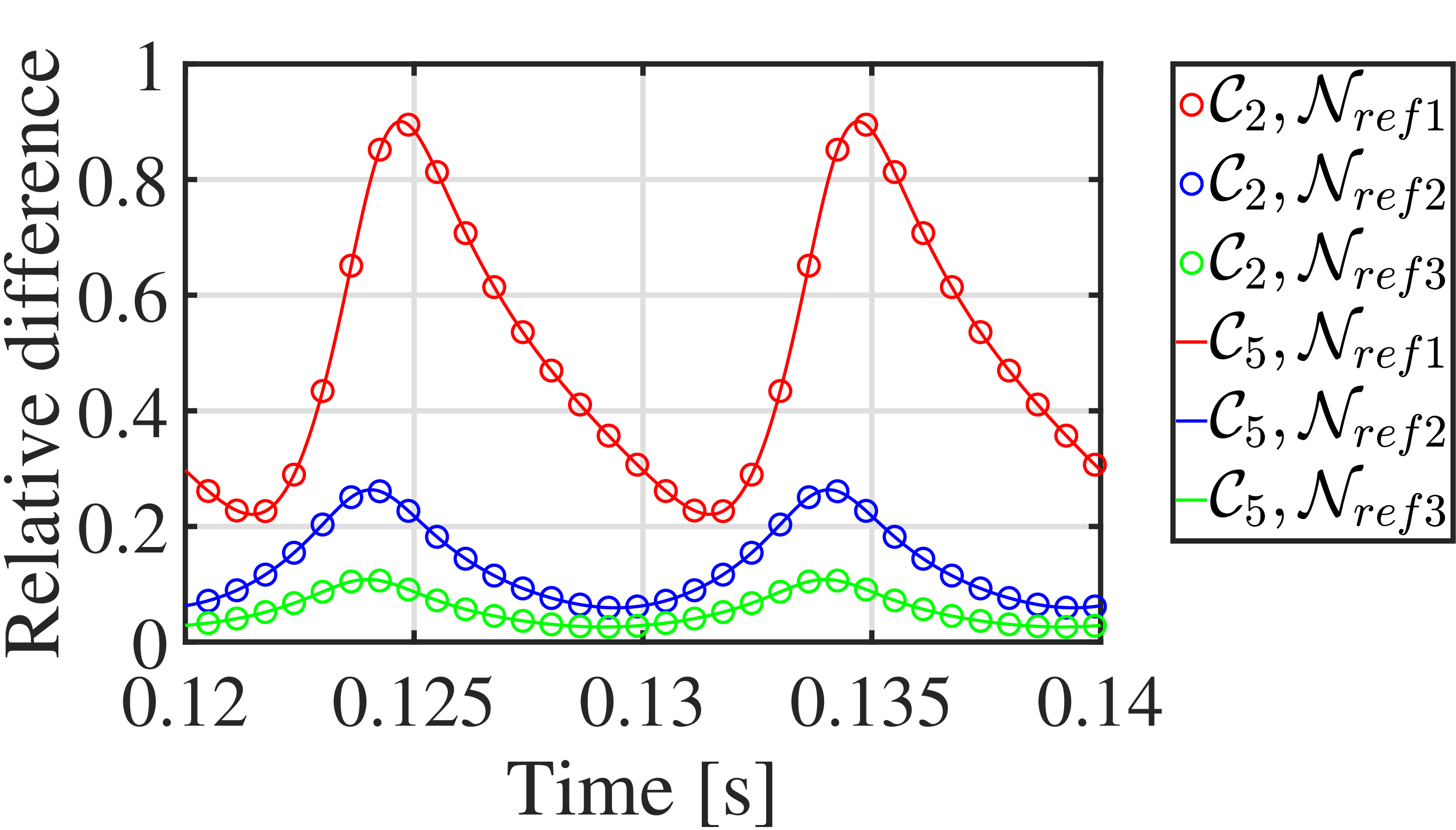}
\label{Busetto:fig:BusettoCurrentDensityIntersection1_NOEC}}
\subfloat[][Copper]{
\includegraphics[width=0.48\textwidth]{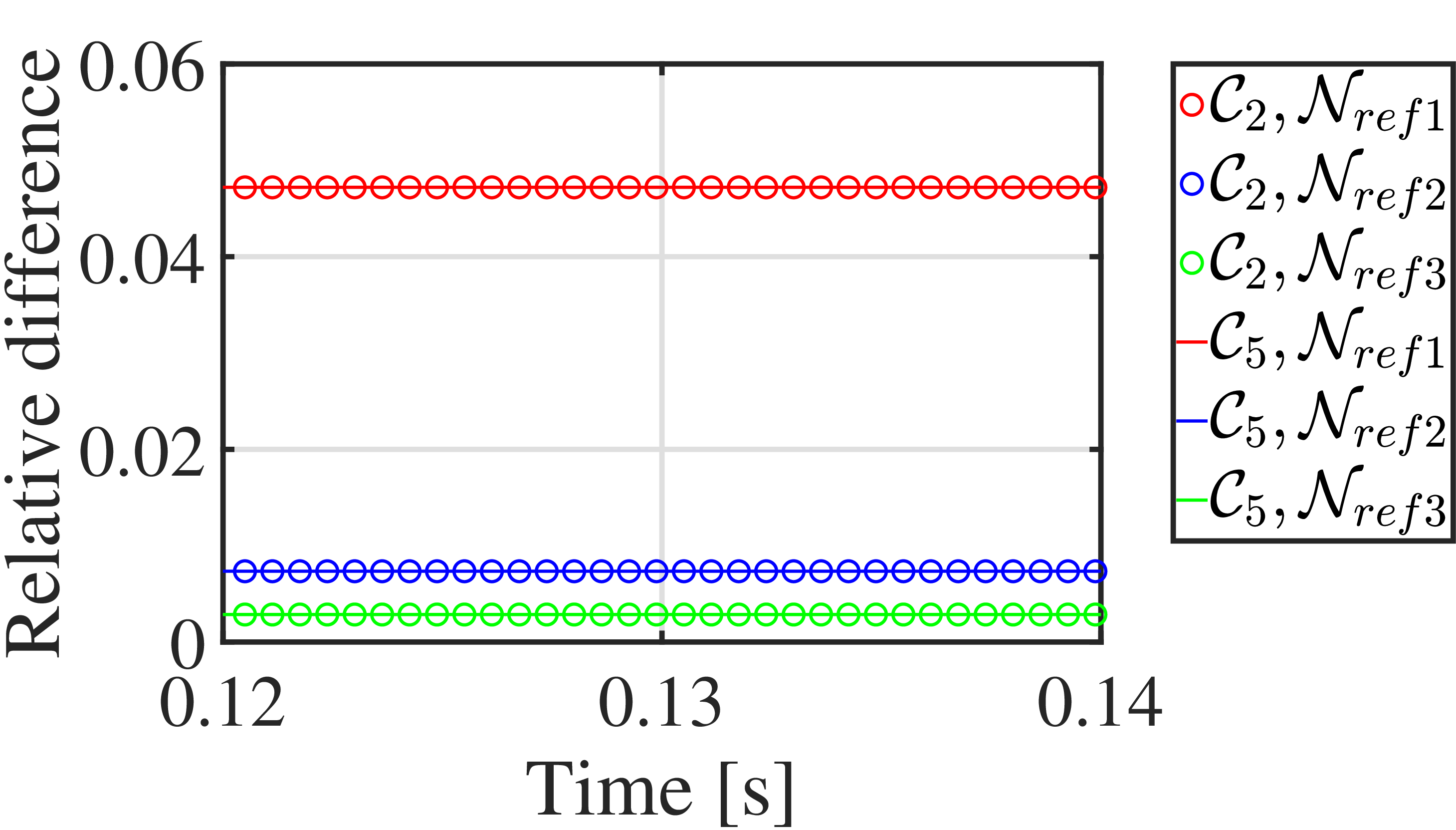}
\label{Busetto:fig:BusettoCurrentDensityIntersection2_NOEC}}
\caption{Test case 3: Evolution in time of $L^2$-norm of difference in current density between numeric solution on three-portion cylinder and analytic solution on infinite cylinder: (a) in first (or third) portion of inner cylinder (iron) and (b) in second portion of inner cylinder (copper).}
\label{Busetto:fig:BusettoCurrentDensityIntersections_NOEC}
\end{figure}

We consider the numerical solution on the three-portion cylinder in two cross sections, one located at the port boundary of one iron region, and the other one in the middle of the copper region, and compare it to the analytic solution of the respective equation on an infinite cylinder of the respective material.
In Figs. \ref{Busetto:fig:BusettoCurrentDensityIntersections} and \ref{Busetto:fig:BusettoCurrentDensityIntersections_NOEC} we report the evolution in time of the $L^2$-norm of the difference between the two solutions for both Test case 2 and Test case 3.
Cylinders $\mathcal{C}_j$, $j\in\{2,5\}$ are considered and all the three refinements $\mathcal{N}_{\mathrm{ref}k}$, $k\in\{1,2,3\}$ for each given cylinder.
We can see that for a given level of refinement the error obtained for the cylinder of different lengths are comparable (the plots overlap).
Moreover, we can see that refining the mesh the differences decrease.
A similar analysis was performed for the magnetic fields showing similar results.

\begin{figure}[htbp]
\centering
\includegraphics[width=0.6\textwidth]{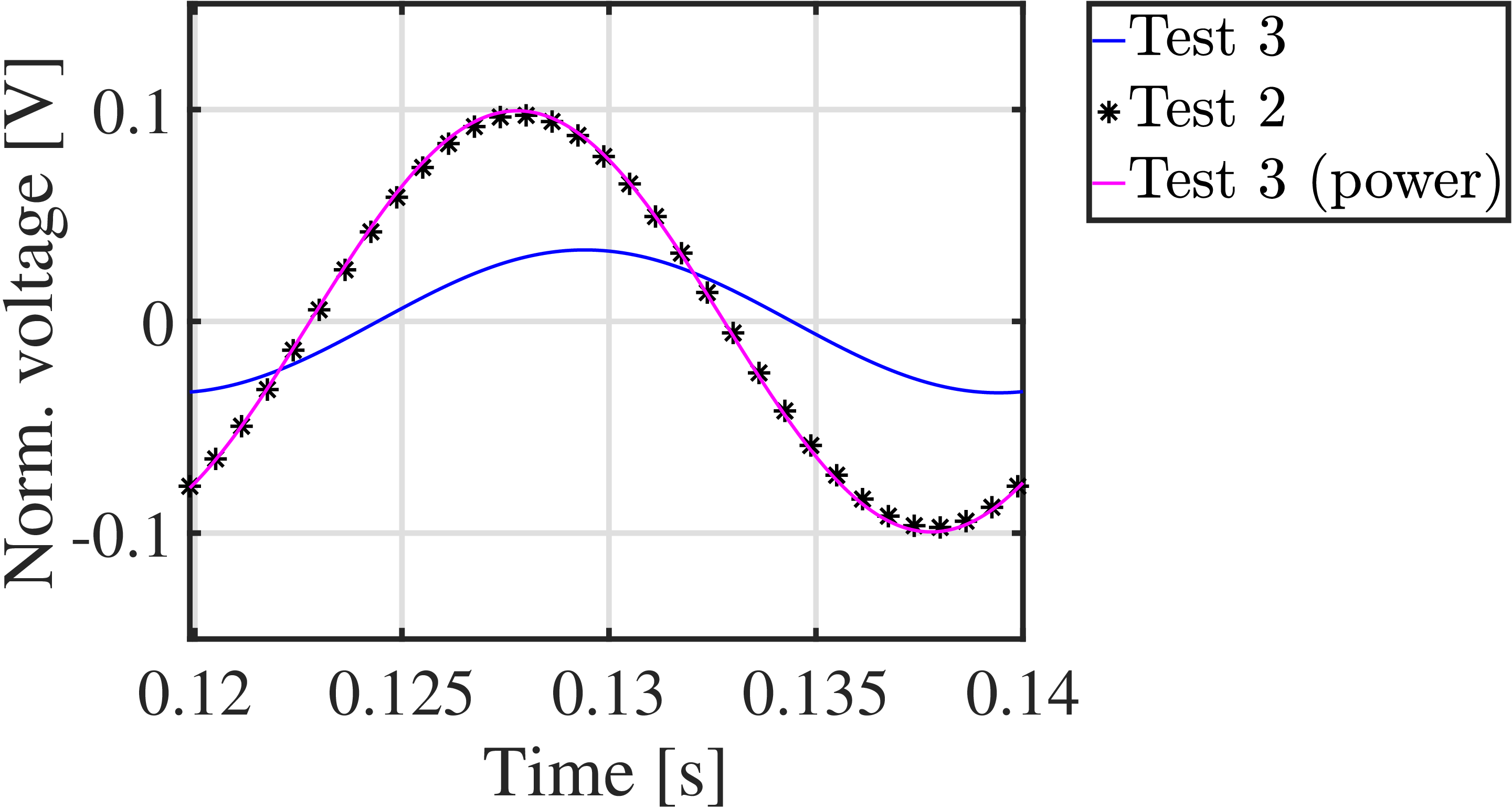}
\caption{Test cases 2 and 3: cylinder $\mathcal{C}_5$ in refinement $\mathcal{N}_{\mathrm{ref}3}$: comparison of the normalized voltage computed using (\ref{Busetto:eq:total_voltage}) for Test cases 2 and 3, and the voltage reconstructed from the power for Test case 2.}
\label{Busetto:fig:Three_portion_80}
\end{figure}

Considering the voltage, we expect the values for Test case 3 to be similar to the values for Test case 2.
Indeed, even when eddy currents are present in the copper part, their effects should be small (see Section \ref{Busetto:sec:ChoiceOfPhysicalParameters}). 
In Fig. \ref{Busetto:fig:Three_portion_80}, we plot the normalized voltage $V_1$ (i.e., the voltage divided by the number of layers of the considered mesh) computed using (\ref{Busetto:eq:total_voltage}) for both Test case 2 and Test case 3 and the voltage reconstructed from the power for Test case 3.
We observe that in case eddy currents are not present everywhere in the inner cylinder, the voltage reconstructed from (\ref{Busetto:eq:total_voltage}) is very different from the one obtained in Test case 2.
Conversely, the voltage reconstructed from the power almost coincides with the one of Test case 2 in agreement with what was explained in Section \ref{Busetto:sec:VoltageReconstruction}.

\begin{figure}[htbp]
\centering
\subfloat[][Test case 2]{
\includegraphics[width=0.48\textwidth]{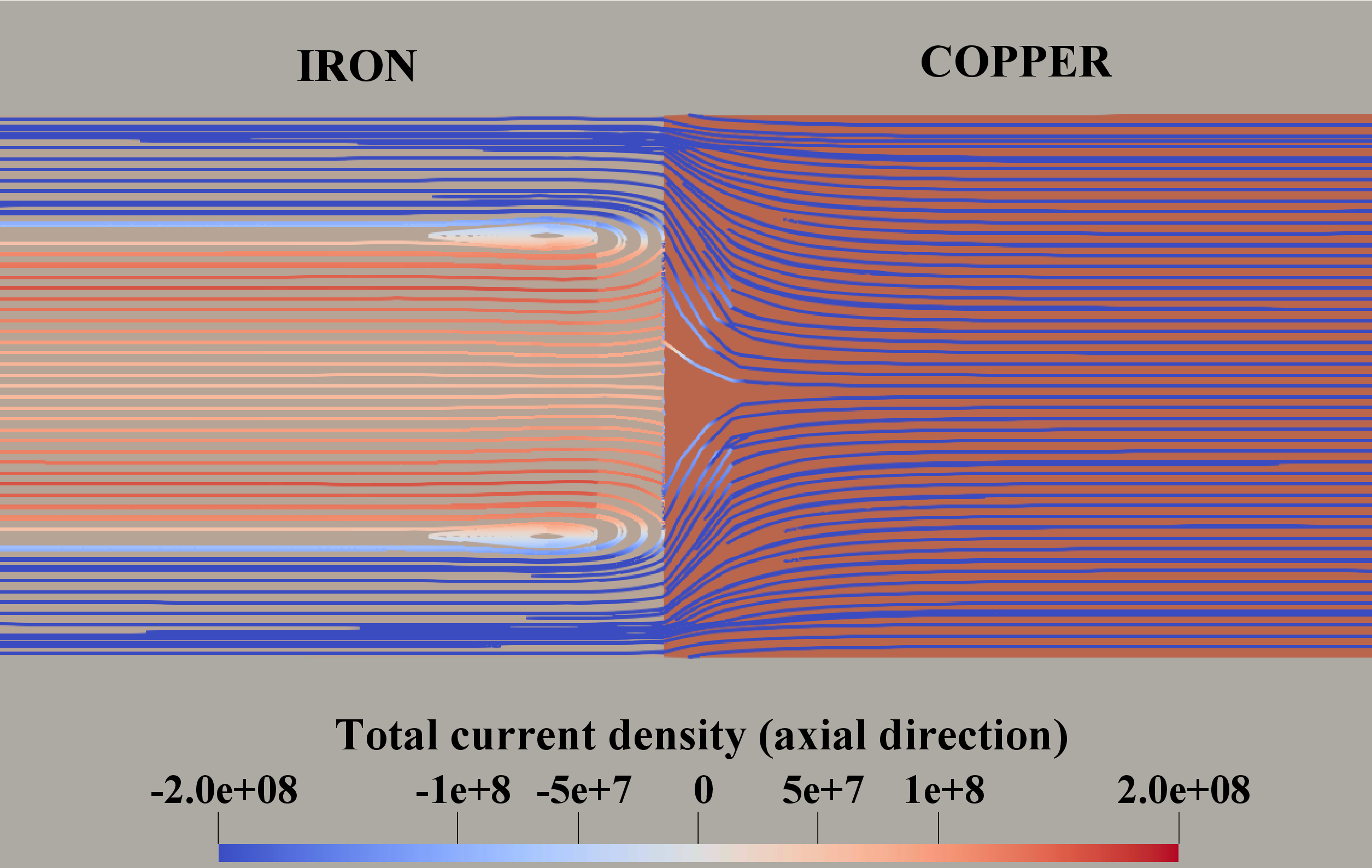}
\label{Busetto:fig:streamlines_EC}}
\subfloat[][Test case 3]{
\includegraphics[width=0.48\textwidth]{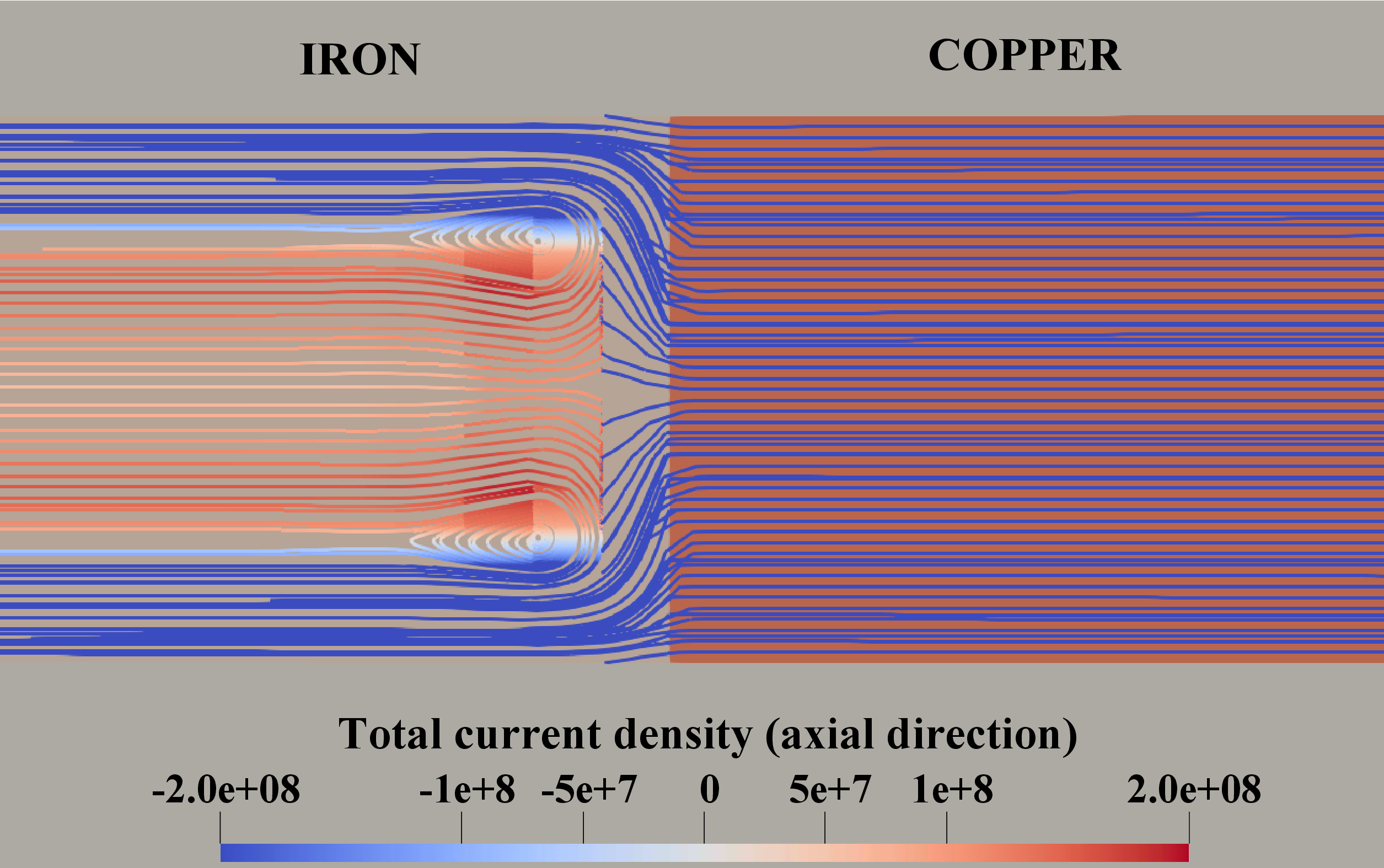}
\label{Busetto:fig:streamlines_NOEC}}
\caption{Current density streamlines: (a) Test case 2, (b) Test case 3.}
\label{Busetto:fig:streamlines}
\end{figure}

In Fig.~\ref{Busetto:fig:streamlines} we compare qualitatively the streamlines of the current density along the axial direction for Test case 2 (Fig.~\ref{Busetto:fig:streamlines_EC}) and for Test case 3 (Fig.~\ref{Busetto:fig:streamlines_NOEC}) at the final time instant $T$.
In particular, we focus on the area near the intersection of the first portion (iron) and the second portion (copper) of the inner cylinder.
As expected, we can see that in the former case the transition happens across the iron and the copper part, whereas in the latter case it happens only in the iron part, as the streamline distribution inside the copper is constant. In both cases we can see the presence of eddies in the iron part. 

\section{Conclusion}\label{Busetto:sec:Conclusion_and_Outlook} 
In this paper we have introduced a two-domain two-step approach for the solution of the eddy current equations in a domain with two ports where the eddy current effects can be neglected in a conductive subdomain.
This approach allows to efficiently compute eddy currents in any situation where their effects are negligible in a subdomain in which their computation would be expensive, e.g., when interpolation would be required due to remeshing of a moving subdomain.
We addressed how the voltage can be reconstructed in this approach.
Numerical results confirm the validity of the proposed method. 

\subsubsection*{Acknowledgement}
This work was supported by ABB Corporate Research.

%
%


\begin{thebibliography}{6}
%
\bibitem{Busetto:bib:ArcSimLV} Bianchetti, R., Adami, A., Fagiano, L., Gati, R., Hofstetter, L.:
Arc Simulation in Low Voltage Switching Devices, a Case Study.
Plasma Physics and Technology, 2, 1, 5-8 (2015).

\bibitem{Busetto:bib:Eddy_current_effects} Chadebec, O., Meunier, G., Mazauric, V., Le Floch, Y., Labie, P.:
Eddy-Current Effects in Circuit Breakers During Arc Displacement Phase.
IEEE Transactions on Magnetics, Institute of Electrical and Electronics Engineers, 40 2, 1358--1361 (2004). \url{doi:10.1109/TMAG.2004.824768}

\bibitem{Busetto:bib:eddy_currents_arc_motion}
Yin, J., Wang, Q., Li, X., Xu, H.:
Numerical Study of Influence of Frequency and Eddy Currents on Arc Motion in Low-Voltage Circuit Breaker.
IEEE Transactions on Components, Packaging and Manufacturing Technology, 8, 8, 1373--1380 (2018). \url{doi:10.1109/TCPMT.2018.2849985}

\bibitem{Busetto:bib:Frequency_stable} Ostrowski, J., Hiptmair, R.:
Frequency-Stable Full Maxwell in Electro-quasistatic Gauge.
SIAM Journal on Scientific Computing, 43, 4, B1008--B1028 (2021). \url{doi:10.1137/20M1356300}

\bibitem{Busetto:bib:Electric_Circuit_Element} Ciuprina, G., Ioan, D., Sabariego, R.V.:
Electric circuit element boundary conditions in the finite element method for full-wave passive electromagnetic devices.
J.~Math.~Industry  12, 7 (2022). \url{doi:10.1186/s13362-022-00122-1}

\bibitem{Busetto:bib:Nodal_auxiliary} Hiptmair, R., Xu, J.:
Nodal auxiliary space preconditioning in H (curl) and H (div) spaces.
SIAM Journal on Numerical Analysis, 45, 6, 2483--2509 (2007). \url{doi:10.1137/060660588}

\bibitem{Busetto:bib:Weeks_Transmission} Weeks, W.L.:
Transmission and Distribution of Electrical Energy.
Harper \& Row, London (1981).

\bibitem{Busetto:bib:Monk} Monk, P.:
Finite Element Methods for Maxwell's Equations.
Clarendon Press, Oxford (2003).

\end{thebibliography}
\end{document}